\documentclass{elsart}
\journal{Theoretical Computer Science}

\usepackage[applemac]{inputenc}
\usepackage{bm,amsfonts,amsmath,amssymb,epsfig,times,amscd,latexsym}
\usepackage{graphicx}
\usepackage{epstopdf}

\def\Dar{\hbox{$\downarrow$}}
\let\e=\varepsilon

\def\gL{\mathbin{\scriptstyle\wedge}}
\def\jR{\mathbin{{\scriptstyle\vee}}}
\def\jjR{\mathbin{\bigvee}}

\def\resp{\hbox{\it resp. }}

\def\hh{\mathbf h}
\def\HH{\mathcal H}
\def\NN{\mathbf N}

\def\xx{{x}}
\def\yy{{y}}
\def\zz{{z}}
\def\ww{{w}}

\def\Irr{\Sigma}
\def\length#1{\hbox{$|#1|$}}

\begin{document}

\begin{frontmatter}



\title{Finite transducers for divisibility monoids}


\author{Matthieu PICANTIN}
\address{LIAFA CNRS UMR 7089\\
Université Paris 7 Denis Diderot\\
F-75251 Paris France}
\ead{picantin@liafa.jussieu.fr}
\ead[url]{http://www.liafa.jussieu.fr/$\sim$picantin}

\begin{abstract}
Divisibility monoids are a natural lattice-theoretical generalization of Mazurkiewicz trace monoids, namely monoids in which the distributivity of the involved divisibility lattices is kept as an hypothesis, but the relations between the generators are not supposed to necessarily be commutations. Here, we show
that every divisibility monoid admits an explicit finite transducer which allows to compute
normal forms in quadratic time. In addition, we prove that every divisibility monoid is
biautomatic.
\end{abstract}

\end{frontmatter}

\section{Introduction}
\label{Sec:Introduction}

\noindent The goal of this paper is to establish the following result~: 

\smallbreak
\noindent
{\bf Main Theorem.} {\sl Every left divisibility monoid admits an explicit finite transducer which
allows to compute right normal forms in quadratic time.}

\smallbreak\noindent Mazurkiewicz's trace theory provides a well-investigated mathematical model
for the sequential behavior of a parallel system in which the order of two independent actions is
regarded as irrelevant. This is achieved by considering a free partially commutative monoid, namely,
the free monoid of all words over a fixed alphabet modulo the congruence generated by equations of
the form~$ab=ba$ for pairs of independent actions~$(a,b)$. Roughly speaking, a letter corresponds
to an event and two letters commute when the corresponding events can occur simultaneously.
However, there are several areas in computer science where one would like to consider more general
equations of the form~$ab = cd$, rather than just~$ab = ba$ as in trace theory.  Left divisibility
monoids have been introduced as a natural algebraic generalization of Mazurkiewicz's trace monoids,
namely monoids in which the distributivity of the underlying left divisibility lattices is kept
as an hypothesis, but the relations between the generators are not supposed to necessarily be
commutations. 

The purpose of this paper is to study how to compute efficiently normal forms in left divisibility monoids. Following Thurston's original idea about the automaticity of the braid groups (see~\cite{eps}), we shall construct an explicit finite transducer---that is a finite automaton with output---allowing to compute normal forms in every left divisibility monoid. Since a standard transducer reads words \emph{from the left to the right}, we need to define a \emph{right} normal form, even if a \emph{left} normal form for
elements in a \emph{left} divisibility monoid---like the normal form defined by Kuske \cite{kus} and
generalizing the Cartier-Foata normal form known from the theory of Mazurkiewicz traces \cite{caf}---seems to be \emph{a priori} the most pertinent choice. 

Before describing the transduction machinary, we exhibit several nice properties of this new normal form with, in particular, a deep geometric property concerning the associated so-called Cayley graph. Our work provides a detailed and complete proof to Kuske's claim that every left divisibility monoid is
automatic~\cite{kus}. Furthermore, we show that every left divisibility monoid is
(both left and right) biautomatic, according to Hoffmann's terminology~\cite{hof}. 

The rest of this paper is organized as follows. In Section~\ref{Sec:Background1}, we recall several definitions about automaticity for monoids. In Section~\ref{Sec:Background2}, we gather the needed basic properties of left divisibility monoids. Section~\ref{Sec:NormalForm} introduces the right normal form. We show that the language of right normal forms has good properties, preparing the proof of how it provides a biautomatic structure to every left divisibility monoid. In Section~\ref{Sec:Transducers}, we then state the main results of this paper (Theorems~\ref{T:Biautomatic} and~\ref{T:Transducer}), discuss and illustrate~them.

\section{Background from automaticity of monoids}
\label{Sec:Background1}

\noindent In this section, we review the theory of automatic monoids---we may focus on cancellative monoids, since left divisibility monoids are defined to be cancellative---and, according to this context, we recall links between automata and transducers. We refer the reader interesting with the notions of automatic structures to~\cite{bau,eps} for automatic groups and to~\cite{cam,hof,hot,hud,sst} for automatic monoids. A general reference about transducers is~\cite{ber}.

\subsection{Languages, automata and transducers}

\noindent We first give a brief introduction to formal language theory (particularly regular languages).

\noindent For a finite set~$X$, let~$X^*$ denote the set of all finite words over the alphabet~$X$, including the empty word~$\varepsilon$.  \noindent For a word~$u$, let~$|u|$ denote the length of~$u$ and let~$u_{[t]}$ (\resp $u^{[t]}$) denote its length~$t$ suffix (\resp prefix) for~$t<|u|$ and $u$ itself for~$t\geq|u|$.

\noindent Some computations on words and on languages can be interpreted as a work of a machine,
which being in a state $p$ and receiving as input a letter $x$, goes into a state~$q$ and possibly outputs a word $w$. Such machines are formalized by the following definitions.

\begin{defn}\label{D:Auto} A (deterministic) \emph{automaton} is a set~${\mathcal A}=(X,Q,q_-,Q_+,\tau)$, where
\begin{enumerate}
\item $X$ is a finite set (the \emph{input alphabet}),
\item $Q$ is a set (the set of \emph{states}),
\item $q_-$ is a fixed element in~$Q$ (the \emph{initial state}),
\item $Q_+$ is a fixed subset of~$Q$ (the set of \emph{accepting states}), and
\item $\tau:X\times Q\rightarrow Q$ is a mapping (the \emph{transition function}).
\end{enumerate}
The map~$\tau$ can be extended to~${\bm\tau}:X^*\times Q\rightarrow Q$  by~${\bm\tau}(\e,q)=q$ and~${\bm\tau}(xu,q)={\bm\tau}(u,\tau(x,q))$ for~$x\in X$, $u\in X^*$ and~$q\in Q$. A word~$w$ over~$X$ is \emph{recognized} by~${\mathcal A}$ if ${\bm\tau}(w,q_-)$ belongs to~$Q_+$. 
\end{defn}

\begin{defn}\label{D:Trans} A (sequential) \emph{transducer} is a set~${\mathcal T}=(X,Y,Q,q_-,Q_+,\tau,\lambda)$, where
\begin{enumerate}
\item $(X,Q,q_-,Q_+,\tau)$ is an automaton,
\item $Y$ is a finite set (the \emph{output alphabet}), and
\item $\lambda:X\times Q\rightarrow Y^*$ is a mapping (the \emph{output function}).
\end{enumerate}
The map~$\lambda$ can be extended to~${\bm\lambda}:X^*\times Q\rightarrow Y^*$ by~${\bm\lambda}(\e,q)=\e$ and~${\bm\lambda}(xu,q)=\lambda(x,q){\bm\lambda}(u,\tau(x,q))$ for~$x\in X$, $u\in X^*$ and~$q\in Q$.
\end{defn}

\noindent An automaton (\resp a transducer) is \emph{finite} if the set~$Q$ is finite. A finite automaton (\resp a finite transducer) can be represented as a labelled directed graph, known as a \emph{Moore diagram}. The vertices of such a graph correspond to the states of the automaton (resp. the transducer), and, for every letter~$x$ of the input alphabet~$X$, an arrow labelled by $x$ (\resp by~$x|\lambda(x,q)$) goes from the state~$q\in Q$ to the state~$\tau(x,q)$. An incoming unlabelled arrow represents the initial state. Accepting states are denoted by double circles. 

\begin{defn} A language is \emph{regular} whenever it is the language of words recognized by some finite automaton. 
\end{defn}

\begin{figure}[thb]
	\begin{center}
	\includegraphics[width=230pt]{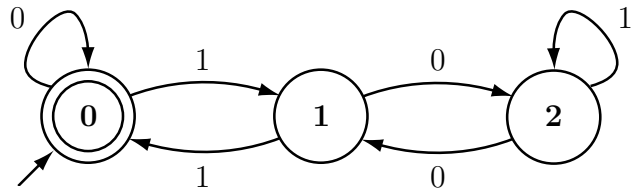}
		\put(0,80){\makebox(0,0){$$}}
		\put(-203,27){\makebox(0,0){$\mathbf{0}$}}
		\put(-115,27){\makebox(0,0){$\mathbf{1}$}}
		\put(-27,27){\makebox(0,0){$\mathbf{2}$}}
		\put(-230,64){\makebox(0,0){$0$}}
		\put(-160,48){\makebox(0,0){$1$}}
		\put(-71,48){\makebox(0,0){$0$}}
		\put(0,64){\makebox(0,0){$1$}}
		\put(-160,4){\makebox(0,0){$1$}}
		\put(-71,4){\makebox(0,0){$0$}}
	\end{center}
\caption{Automaton deciding division by~$3$ in base~2.}
\label{F:Div23}
\end{figure}

\begin{figure}[thb]
	\begin{center}
	\includegraphics[width=230pt]{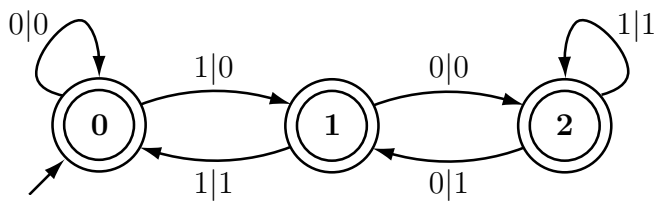}
		\put(0,80){\makebox(0,0){$$}}
		\put(-203,27){\makebox(0,0){$\mathbf{0}$}}
		\put(-115,27){\makebox(0,0){$\mathbf{1}$}}
		\put(-27,27){\makebox(0,0){$\mathbf{2}$}}
		\put(-230,64){\makebox(0,0){$0|0$}}
		\put(-160,48){\makebox(0,0){$1|0$}}
		\put(-71,48){\makebox(0,0){$0|0$}}
		\put(0,64){\makebox(0,0){$1|1$}}
		\put(-160,4){\makebox(0,0){$1|1$}}
		\put(-71,4){\makebox(0,0){$0|1$}}
	\end{center}
\caption{Transducer allowing conversion from base~2 to base~3.}
\label{F:Base23}
\end{figure}

\begin{exmp}\label{E:Base} For given positive integers~$b$ and~$k$, divisibility by~$k$ in base~$b$ can be decided by a finite automaton, and it turns out that converting integers from base~$b$ to base~$k$ can be made by using a finite transducer~: the latter computes the remainder---the final state---and the quotient---the output---modulo~$b$. Figure~\ref{F:Div23} displays an automaton which reads from the left to the right and decides whether a base~2 integer is divisible by~3. Figure~\ref{F:Base23} displays the associated transducer which reads from the left to the right and allows---via multiple runs---to convert an integer from base~2 to base~3.
\end{exmp}

\begin{rem} Most of the transducers we shall consider have only accepting states. This special feature allows to use these transducers iteratively.
\end{rem}
\subsection{Notions of automatic monoids}

\noindent First defined by Thurston two decades ago, automatic groups attracted a lot of attention in geometric and combinatorial group theory and are the subject of a major book~\cite{eps} (see also~\cite{bau}). Roughly speaking, an automatic group is a finitely generated group for which one can check, by means of a finite automaton, whether two words over a finite generating alphabet represent the same element or not, and whether or not the elements they represent differ by multiplication by a single generator. A few years ago, the notion of automaticity was generalized for semigroups and monoids~: it is worth mentioning the work by~Hoffmann in~\cite{hof} (see also~\cite{cam,hot,hud,otr}).

As with automatic groups, we may consider automata reading pairs of words where we introduce a \emph{padding symbol} to deal with the case where the lengths of the two words are not the same. One can introduce the paddings on the right or on the left.

\begin{defn} For every alphabet~$X$, the mappings~$\overrightarrow{\phantom{;}.\phantom{;}}_X$ and~$\overleftarrow{\phantom{;}.\phantom{;}}_X$ from~$X^*\times X^*$ to~$Y^*$ with~$\$\not\in X$ and~$Y=(X\cup\{\$\})\times(X\cup\{\$\})\setminus\{(\$,\$)\}$ are defined by\[\overrightarrow{(x_1\cdots x_n,y_1\cdots y_m)}_X=\left\{\begin{array}{ll}
(x_1,y_1)\cdots(x_n,y_n)&\hbox{for~}n=m,\\
(x_1,y_1)\cdots(x_n,y_n)(\$,y_{n+1})\cdots(\$,y_m)&\hbox{for~}n<m,\\
(x_1,y_1)\cdots(x_m,y_m)(x_{m+1},\$)\cdots(x_n,\$)&\hbox{for~}n>m,
\end{array}
\right.
\]and~$\overleftarrow{(x_1\cdots x_n,y_1\cdots y_m)}_X$ to be the mirror of~$\overrightarrow{(x_n\cdots x_1,y_m\cdots y_1)}_X$,
where the~$x_i$'s and the~$y_j$'s belong to~$X$ for~$1\leq i\leq n$ and~$1\leq j\leq m$.
\end{defn}

\noindent Hoffmann purposed then four notions of automaticity for semigroups~: roughly speaking, for~$\pi,\mu$ in~$\{{\rm left},{\rm right}\}$, a semigroup is said to be $\pi$-$\mu$-automatic if it is automatic with $\pi$ the direction of \emph{padding} and $\mu$ the direction of \emph{multiplication}.

\noindent For a monoid~$M$ generating by a set~$X$, there is a canonical mapping~$\overline{\phantom{;}.\phantom{;}}:X^*\twoheadrightarrow M$.

\begin{defn} Assume that~$M$ is a monoid---or a semigroup---generating by a finite set~$X$ and that~$L$ is a language over~$X$ that maps onto~$M$. Then

$(X,L)$ is a \emph{left-left automatic structure} for~$M$\\
{\phantom{~~~~~~~~~~~~~~~~~~~~~}if~${~}_x^\$L=\{\overleftarrow{(u,v)}_X:u,v\in~L, x\overline{u}=\overline{v}\}$ is regular for~$x\in X\cup\{\varepsilon\}$~;}

$(X,L)$ is a \emph{right-left automatic structure} for~$M$\\
{\phantom{~~~~~~~~~~~~~~~~~~~~~}if~${}_xL^\$=\{\overrightarrow{(u,v)}_X:u,v\in~L, x\overline{u}=\overline{v}\}$ is regular for~$x\in X\cup\{\varepsilon\}$~;}

$(X,L)$ is a \emph{left-right automatic structure} for~$M$\\
{\phantom{~~~~~~~~~~~~~~~~~~~~~}if~${~}^\$L_x=\{\overleftarrow{(u,v)}_X:u,v\in~L, \overline{u}x=\overline{v}\}$ is regular for~$x\in X\cup\{\varepsilon\}$~;}

$(X,L)$ is a \emph{right-right automatic structure} for~$M$\\
{\phantom{~~~~~~~~~~~~~~~~~~~~~}if~$L_x^\$=\{\overrightarrow{(u,v)}_X:u,v\in~L, \overline{u}x=\overline{v}\}$ is regular for~$x\in X\cup\{\varepsilon\}$.}
Those automata accepting such languages are called \emph{equality recognizer automata} for~$x=\varepsilon$ and~\emph{multiplier automata} for~$x\in X$.
\end{defn} 

\noindent The notion of \emph{automatic} as defined in~\cite{cam} for semigroups is equivalent to the notion of \emph{right-right automatic} here. These four notions of automaticity are shown to be independent for general semigroups and to collapse into a dual notion of $\mu$-automaticity for cancellative monoids (whether automaticity implies biautomaticity is still an open question for groups). Roughly speaking, the property of automaticity for a cancellative monoid does not depend on the direction of padding or reading.

\begin{prop}\label{P:RL} Assume that~$M$ is a cancellative monoid generating by a finite set~$X$ and~$L$ is a language over~$X$ that maps onto~$M$. Then, for~$\mu$ in~$\{{\rm left},{\rm right}\}$, $(X,L)$ is a left-$\mu$ automatic structure for~$M$ if and only if $(X,L)$ is a right-$\mu$ automatic structure for~$M$.
\end{prop}

\begin{rem} Hoffmann exhibited an example of a cancellative monoid that satisfies  all four of notions of automaticity but which is not right-biautomatic, that is, which does not admit a structure being both right-left and right-right automatic.
\end{rem}

\subsection{Hoffmann's criterium}

\noindent The automatic structures for groups are characterized by a geometric condition on the associated Cayley graph, known as the \emph{fellow traveller property} (see~\cite[Theorem 2.3.5]{eps}). Roughly speaking, a formal language~$L$ over an alphabet~$X$ generating a group~$G$ satisfies the fellow traveller property if paths in the Cayley graph of~$G$ with respect to~$X$, which are labelled by words in~$L$ and which eventually converge to within a distance of~$1$, never diverge beyond some given distance. If one seeks to apply this condition to monoids, one must decide what one means by \emph{distance} in a monoid Cayley graph. For cancellative monoids, a convenient notion is the following~:

\begin{defn} Assume that~$M$ is a monoid with a finite generating alphabet~$X$. Then the function~$d_X:M\times M\rightarrow{\mathbb N}\cup\{\infty\}$ (\resp $d^X$) defined by
\[d_X(a,b)=\min\{|w|:w\in X^*\hbox{\rm ~and~} (\overline{w}a=b \hbox{\rm ~or~} a=\overline{w}b)\}\]
(\resp $d^X(a,b)=\min\{|w|:w\in X^*\hbox{\rm ~and~} (a\overline{w}=b \hbox{\rm ~or~} a=b\overline{w})\}$) is called the \emph{left (\resp right) directed distance function of~$M$ with respect to~$X$}.
\end{defn}

\noindent Note that such a notion of distance does not satisfy the triangular inequality.

\begin{defn} Assume that~$M$ is a monoid with a finite generating alphabet~$X$. Then a language~$L$ over~$X$ is said to satisfy the \emph{left (\resp right) directed fellow traveller property (with respect to~$M$)} if there exists a positive integer~$k$ such that, for any two words~$u,v$ in~$L$ satisfying~$d_X(u,v)\leq 1$ (\resp satisfying~$d^X(u,v)\leq 1$), we have~$d_X(u_{[t]},v_{[t]})<k$ (\resp $d^X(u^{[t]},v^{[t]})<k$) for every nonnegative integer~$t$. 
\end{defn}

\noindent What we refer to Hoffmann's criterium is the following result. Let us mention that the original version~\cite[Proposition 8.3]{hof} is stated in term of a semigroup~$S$ such that, for every~$a,b,d$ in~$S$ satisfying~$ab=ad$, $cb=cd$ holds for every~$c$ in~$S$. Now, in the case of a monoid, the previous hypothesis is equivalent to left cancellativity. 

\begin{thm}\label{P:HC} Assume that~$M$ is a right (\resp left) cancellative monoid with a finite generating alphabet~$X$. Then every regular language~$L$ over~$X$ mapping onto~$M$ and satisfying the left (\resp right) directed fellow traveller property with respect to~$M$ provides a left (\resp right) automatic structure for~$M$.
\end{thm}

\noindent Let us mention that several different geometric conditions characterizing automatic monoids were investigated (see~\cite{sst} for instance).

\section{Background from left divisibility monoids}
\label{Sec:Background2}

\noindent In this section, we list some basic properties of left divisibility monoids, and summarize results by Droste \&~Kuske about them. For all the results quoted in this section, we refer the reader to~\cite{dku,kur,kus}.

\subsection{Divisors and multiples in a monoid}
Assume that~$M$ is a monoid. We say that~$M$ is \emph{conical} if~$1$ is the only invertible element
in~$M$. For~$a,b$ in~$M$, we say that~$b$ is a left divisor of~$a$---or that~$a$ is a right
multiple of~$b$---if~$a=bd$ holds for some~$d$ in~$M$. The set of the left divisors of~$b$ is
denoted by~$\Dar(b)$. An element~$c$ is a right lower common multiple---or a right lcm---of~$a$
and~$b$ if it is a right multiple of both~$a$ and~$b$, and every right common multiple of~$a$
and~$b$ is a right multiple of~$c$. Right divisor and left multiple are defined symmetrically. 

\noindent If~$c$, $c'$ are two right lcm's of~$a$ and~$b$, necessarily~$c$ is a left
divisor of~$c'$, and~$c'$ is a left divisor of~$c$. If we assume~$M$ to be conical and
cancellative, we have~$c=c'$. In this case, the unique right lcm of~$a$ and~$b$ is denoted by~$a\jR
b$. Cancellativity and conicity imply that left and right divisibility are order relations.

\noindent Let~$M$ be a monoid. An \emph{irreducible element} of~$M$ is defined to be a non
trivial element~$a$ such that~$a=bc$ implies~$b=1$ or~$c=1$. The set of the irreducible elements
in~$M$ can be written as~$(M\setminus\{1\})\setminus(M\setminus\{1\})^2$.

\subsection{Main definitions and properties for left divisibility monoids}

\noindent Let~$(P,\leq)$ be a partially ordered set. Then, for any~$a$ in~$P$, $\Dar(a)$ comprises
all elements dominated by~$a$, that is, $\Dar(a)=\{b\in P;b\leq a\}$. The \emph{width} of a partially ordered set is the maximal size of an antichain, that is, a subset such that any two distinct elements are incomparable. The partially ordered set
$(P,\leq)$ is a lattice if, for any two~$a,b$ in~$P$, the least upper bound~$\sup(a,b)=a\jR b$ and
the largest lower bound~$\inf(a,b)=a\gL b$ exist. The lattice $(P,\leq)$ is \emph{distributive} if
$a\gL(b\jR c)=(a\gL b)\jR(a\gL c)$ for any~$a,b,c$ in~$P$. This is equivalent to $a\jR(b\gL
c)=(a\jR b)\gL(a\jR c)$ for any~$a,b,c$ in~$P$. For properties of finite distributive lattices, we
refer the reader to~\cite{brk}.

\begin{defn}\label{D:DiviMonoid} A monoid~$M$ is called a \emph{left divisibility monoid}---or simply
a \emph{divisibility monoid}---if~$M$ is cancellative and finitely generated by its irreducible
elements, if any two elements admit a left gcd and if~every element~$a$ dominates a
finite distributive lattice~$\Dar(a)$.
\end{defn}

\noindent Note that the finiteness requirement on lattices is in fact not necessary since it follows from the
other stipulations. Note also that cancellativity and the lattice condition imply conicity. The left gcd of two
elements~$a,b$ will be denote by~$a\gL b$. The \emph{length}~$\length{a}$ of an
element~$a$ is defined to be the height of the lattice~$\Dar(a)$.

\begin{exmp}\label{E:DiviMonoid} Every (finitely generated) trace monoid is a divisibility
monoid. Both the monoids~$\langle~\xx,\yy,\zz:\xx\yy=\yy\zz~\rangle$
and~$\langle~\xx,\yy,\zz:\xx^2=\yy\zz,~\yy\xx=\zz^2~\rangle$ are not trace but
left divisibility monoids. The
monoid~$\langle~\xx,\yy,\zz:\xx^2=\yy\zz,~\xx\yy=\zz^2\rangle$ is neither a left nor
a right divisibility monoid---a monoid being called a \emph{right divisibility monoid} if its antiautomorphic image is a left divisibility monoid.
\end{exmp}

\noindent An easy but crucial fact about left divisibility monoids is the
following.

\begin{lem}\label{L:Lcm} Assume that~$M$ is a left divisibility monoid. Then finitely
many elements in~$M$ admitting at least a right common multiple admit a unique right lcm.
\end{lem}

\noindent The following result states that there exists a decidable class of
presentations that gives rise precisely to all left divisibility monoids.

\begin{thm}\label{T:Kuske}
Assume that~$M$ is a monoid finitely generated by the set~$\Irr$ of its irreducible
elements. Then~$M$ is a left divisibility monoid if and only if 

\noindent(i) $\Dar(xyz)$ is a distributive lattice, 

\noindent(ii) $xyz=xy'z'$ or~$yzx=y'z'x$ implies~$yz=y'z'$, 

\noindent(iii) $xy=x'y'$, $xz=x'z'$ and~$y\not=z$ imply~$x=x'$,

\noindent for any~$x,y,z,x',y',z'$ in~$\Irr$, and if

\noindent(iv) we have~$M\cong\Irr^*\!/\!\!\sim$, with~$\sim$ the congruence
on~$\Irr^*$ generated by the pairs~$(xy,zt)$ for~$x,y,z,t$ in~$\Irr$ and $xy=zt$.
\end{thm}

\noindent Kuske studied a left normal form generalizing the Cartier-Foata normal form known from the theory of Mazurkiewicz traces. This left normal form can be computed by an infinite transducer and Kuske claims in~\cite{kur} that the latter would allow to prove that every left divisibility monoid is automatic. His main result is that the transducer is finite if and only if the monoid is width-bounded, if and only if the monoid is a regular monoid~\cite{sak}. We shall come back to Kuske's \emph{infinite} transducer in Remark~\ref{R:Kuske}.

\section{A right normal form}
\label{Sec:NormalForm}

\noindent Our goal being to construct finite transducers allowing to compute normal forms and standard
transducers reading words from the left to the right, we shall define a \emph{right} normal form. The right
normal form of an element will be defined as a unique decomposition into a
product of so-called hypercubes, where the rightmost hypercube is the maximal
one, and so on.

\noindent Our aim is to show that the right normal form we consider is associated with a biautomatic structure using Hoffmann's criterium.

\subsection{Definition of a right normal form}

\noindent A natural \emph{left} normal form for \emph{left} divisibility monoids is defined in~\cite{kus} and then called \emph{Foata normal form}. The latter does not seem to be always the best fitted to standard transducers, which reads from the left to the right. One could work with \emph{right} divisibility monoids, but this would devalue the property for left-right reading of being standard. A convenient choice is to construct a \emph{right} normal form for \emph{left} divisibility monoids. Although less natural \emph{a priori}, this choice will be shown to provide equivalent features.

\begin{defn} Assume that~$M$ is a left divisibility monoid. An element~$h$ in~$M$ is called a \emph{hypercube} if there exist irreducibles~$x_1,\ldots,x_n$ satisfying~$h=x_1\jR\cdots\jR x_n$. By
convention, the trivial element~$1$ is a hypercube.
\end{defn}

\noindent Since every finite distributive lattice whose upper bound is the join of its atoms is a
hypercube (see~\cite{brk} or for instance~\cite[page 107]{sta}), a hypercube in a left
divisibility monoid is an element~$h$ whose lattice~$\Dar(h)$ is a hypercube. In particular, since every interval of a hypercube lattice is a hypercube lattice, every divisor of a hypercube in a left divisibility monoid is a hypercube.

\begin{lem}\label{L:RightDividing} Assume that~$M$ is a left divisibility
monoid. Then every element in~$M$ is right-divided by a unique maximal hypercube.\end{lem}
 
\begin{pf} Let~$d$ be an element in~$M$ and~$\{x_1,\ldots,x_p\}$ be the set of those
irreducible elements that divide~$d$ on the right. Since~$d$ is a common left multiple of
the~$x_i$'s, there exists at least one minimal common left multiple of the~$x_i$'s. Such an element
is therefore a hypercube, namely a $p$-cube. Assume now that~$b$ and~$c$ are two distinct $p$-cubes dividing~$d$ on the right. Then the elements~$b',c'$ in~$M$ satisfying~$d=b'b=c'c$ do not admit
a unique right lcm in the lattice~$\Dar(d)$, contradicting Lemma~\ref{L:Lcm}.\qed
\end{pf}

\begin{defn}\label{D:RightNF} Assume that~$M$ is left divisibility monoid. The \emph{right normal form} of a non-trivial element~$a$ in~$M$ is the unique decomposition into non-trivial hypercubes
$\NN(a)=h_p\cdot\ldots\cdot h_1$ such that~$a=h_p\cdots h_1$ holds in~$M$ and~$h_i$ is the unique
maximal hypercube right-dividing~$h_p\cdots h_i$ for~$1\leq i\leq p$. Moreover, we set~$\NN(1)=1$.
\end{defn}

\noindent The remainder of this section is devoted to prove several properties of this right normal form. The latter will allow to finally establish that the language of these right normal forms provides a biautomatic structure to every left divisibility monoid.

\subsection{Regularity of the language of right normal forms}

\noindent The first of the two key points is that the normality of a word is characterized by a local condition, what is captured by the following lemma.

\begin{lem}\label{L:KeyLemma} Assume that~$M$ is a left divisibility monoid
and~$\HH$ is the set of its hypercubes. Let~$\hh:M\rightarrow\HH$ map an
element~$a$ to the maximal hypercube right-dividing~$a$. Then $\hh(ab)=\hh(\hh (a)b)$ holds for
any two elements~$a,b$ in~$M$.
\end{lem}

\begin{pf} We use an induction on the length~$\length b$ of~$b$. For~$\length b=0$, the formula follows
from~$\hh^2=\hh$. Assume~$\length b=1$. Then $b$ is an irreducible element, say~$b=x$.
Let~$y_1,\ldots,y_p$ be the distinct irreducible elements right-dividing~$a$, so
right-dividing~$\hh(a)$ by definition. The distributivity condition implies that, for every~$j$,
there exists at most one irreducible element~$z_{i_j}$ satisfying~$y_jx=t_{i_j}z_{i_j}$ for some
irreducible~$t_{i_j}\not=y_j$. Therefore $\hh(ax)$ and~$\hh(\hh(a)x)$ are the $(q+1)$-cube
right-divided by~$x,y_1,\ldots,y_q$ with~$q\leq p$.  We obtain~$\hh(ax)=\hh(\hh(a)x)$ for every
element~$a$ and every irreducible~$x$ in~$M$.

Assume now~$\length b>1$. Then there exist an irreducible~$x$ and an element~$b'$ in~$M$
satisfying~$b=xb'$, and, by induction hypothesis, we obtain
\[\hh(ab)=\hh(axb')=\hh(\hh(ax)b')=\hh(\hh(\hh(a)x)b')=\hh(\hh(a)xb')=\hh(\hh(a)b),\]
which concludes the induction.\qed\end{pf}
 
\begin{prop}\label{P:LocalCondition} Assume that~$M$ is a left divisibility monoid and~${\mathcal
H}$ is the set of its hypercubes. Let~$h_1,\ldots,h_p$ belong to~$\HH$. Then
$h_p\cdot\ldots\cdot h_1$ is a right normal form if and only if so is~$h_{i+1}\cdot h_i$ for~$1\leq
i<p$.
\end{prop}

\begin{pf} Using the formula~$\hh(ab)=\hh(\hh(a)b)$ from Lemma~\ref{L:KeyLemma}, we
find~$\hh(h_p\cdots h_i)=\hh(h_{i+1}h_i)$ for~$1\leq i<p$.\qed\end{pf}

\begin{cor}\label{P:Regularity} Assume that~$M$ is a left divisibility monoid.
Then the language of its right normal forms is regular.
\end{cor}

\begin{pf} Our language is over the finite alphabet~$\HH$ of hypercubes in~$M$. It suffices to take the automaton with~$\HH$ as set of states and with a transition from~$a$ to~$b$ whenever~$a\cdot b$ is normal, that is, whenever $\hh(ab)=b$ holds.\qed\end{pf}

\begin{figure}[htb]
	\begin{center}
	\includegraphics[width=4cm]{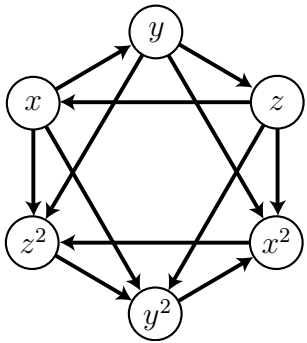}
		\put(0,141){\makebox(0,0){$$}}
		\put(-56.5,117){\makebox(0,0){$y$}}
		\put(-56,10){\makebox(0,0){$y^2$}}
		\put(-103,38){\makebox(0,0){$z^2$}}
		\put(-11,38){\makebox(0,0){$x^2$}}
		\put(-103,90){\makebox(0,0){$x$}}
		\put(-11,90){\makebox(0,0){$z$}}
	\end{center}
\caption{The graph of hypercubes for~$\langle~x,y,z:x^2=yz,y^2=zx,z^2=xy~\rangle$.}
\label{F:Graph}
\end{figure}

\begin{rem}\label{R:Connectivity} Contrary to the case of trace monoids (see \cite[Lemma
3.2]{kro}) or braid monoids (see \cite[Proposition 4.9]{bes}), the \emph{graph of hypercubes} of a general left divisibility monoid---defined
to be the oriented graph with~$\HH$ as set of vertices and with an edge from~$a$ to~$b$
whenever $\hh(ab)=b$ holds---need not necessarily be strongly connected, even if a condition of \emph{irreducibility} (see \cite{pig}) is required. For instance, in the divisibility
monoid $\langle~x,y,z:x^2=yz,y^2=zx,z^2=xy~\rangle$, there is no right normal form like~$z^2\cdot\ldots\cdot x$. The graph of its hypercubes is displayed on
Figure~\ref{F:Graph} (we have omitted the vertices of the two central hypercubes which anyway are
not involved in the strong connectivity). Such a graph is known as the graph of cliques
in~\cite{kro} and as Charney's graph in~\cite{bes}, and its strong connectivity plays a
pivotal r\^ole in the study of trace monoids and braid monoids, respectively.
\end{rem}

\subsection{A fellow traveller property for the language of right normal forms}

\noindent The second key point is that the right normal forms
have a nice behaviour towards both left \emph{and} right multiplication. 

\begin{prop}\label{P:Multiplication} Assume that~$M$ is a left divisibility monoid. Let~$a$
be an element in~$M$ with right normal form~$h_m\cdot\ldots\cdot h_1$ and let~$y,z$ be
hypercubes. Then\begin{itemize}
\item[(i)] the right normal form of~$ya$ is~$h'_m\cdot\ldots\cdot h'_1\cdot y_0$ (or possibly $h'_{m-1}\cdot\ldots\cdot h'_1\cdot y_0$)
with~$y_m=y$,
$y_{i-1}=\hh(y_ih_i)$ and~$y_ih_i=h'_iy_{i-1}$ for~$1\leq i\leq m$.
\item[(ii)] the right normal form of~$az$ is~$z_m\cdot h''_m\cdot\ldots\cdot h''_1$ (or possibly $h''_m\cdot\ldots\cdot h''_1$) with~$z_0=z$,
$h''_i=\hh(h_iz_{i-1})$ and~$z_ih''_i=h_iz_{i-1}$ for~$1\leq i\leq m$.
\end{itemize}
\end{prop}

\begin{figure}[thb]
	\begin{center}
		\includegraphics[width=11cm]{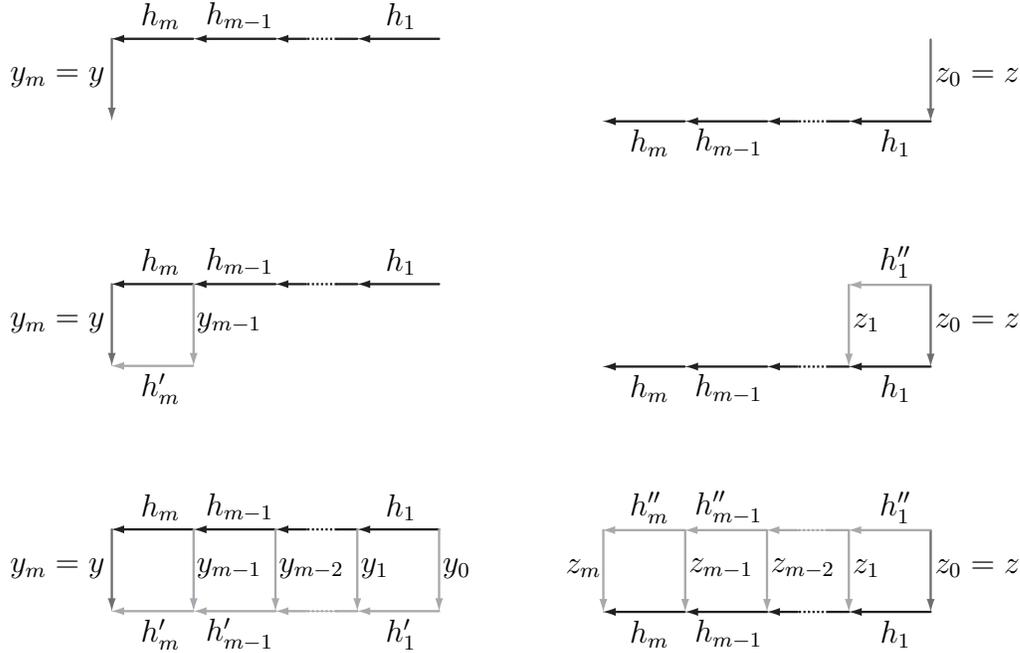}
		\put(0,245){\makebox(0,0){$$}}
		\put(-293,227){\makebox(0,0){$h_m$}}
		\put(-263,227){\makebox(0,0){$h_{m-1}$}}
		\put(-202,227){\makebox(0,0){$h_1$}}
		\put(-332,204){\makebox(0,0){$y_m=y$}}
		\put(16,204){\makebox(0,0){$z_0=z$}}
		\put(-108,179){\makebox(0,0){$h_m$}}
		\put(-78,179){\makebox(0,0){$h_{m-1}$}}
		\put(-15,179){\makebox(0,0){$h_1$}}
		\put(-293,134){\makebox(0,0){$h_m$}}
		\put(-263,134){\makebox(0,0){$h_{m-1}$}}
		\put(-202,134){\makebox(0,0){$h_1$}}
		\put(-15,134){\makebox(0,0){$h''_1$}}
		\put(-332,111){\makebox(0,0){$y_m=y$}}
		\put(-267,111){\makebox(0,0){$y_{m-1}$}}
		\put(-26,111){\makebox(0,0){$z_1$}}
		\put(16,111){\makebox(0,0){$z_0=z$}}
		\put(-293,86){\makebox(0,0){$h'_m$}}
		\put(-108,86){\makebox(0,0){$h_m$}}
		\put(-78,86){\makebox(0,0){$h_{m-1}$}}
		\put(-15,86){\makebox(0,0){$h_1$}}
		\put(-293,41){\makebox(0,0){$h_m$}}
		\put(-263,41){\makebox(0,0){$h_{m-1}$}}
		\put(-202,41){\makebox(0,0){$h_1$}}
		\put(-108,41){\makebox(0,0){$h''_m$}}
		\put(-78,41){\makebox(0,0){$h''_{m-1}$}}
		\put(-15,41){\makebox(0,0){$h''_1$}}
		\put(-332,18){\makebox(0,0){$y_m=y$}}
		\put(-267,18){\makebox(0,0){$y_{m-1}$}}
		\put(-236,18){\makebox(0,0){$y_{m-2}$}}
		\put(-212,18){\makebox(0,0){$y_1$}}
		\put(-181,18){\makebox(0,0){$y_0$}}
		\put(-133,18){\makebox(0,0){$z_m$}}
		\put(-81,18){\makebox(0,0){$z_{m-1}$}}
		\put(-50,18){\makebox(0,0){$z_{m-2}$}}
		\put(-26,18){\makebox(0,0){$z_1$}}
		\put(16,18){\makebox(0,0){$z_0=z$}}
		\put(-293,-7){\makebox(0,0){$h'_m$}}
		\put(-263,-7){\makebox(0,0){$h'_{m-1}$}}
		\put(-202,-7){\makebox(0,0){$h'_1$}}
		\put(-108,-7){\makebox(0,0){$h_m$}}
		\put(-78,-7){\makebox(0,0){$h_{m-1}$}}
		\put(-15,-7){\makebox(0,0){$h_1$}}
	\end{center}
	\vglue 2mm
\caption{The right normal forms of~$ya$ and~$az$ from the right normal form~$h_m\cdot\ldots\cdot h_1$ of~$a$.}
\label{F:Multiplication}
\end{figure}

\begin{pf} (Figure~\ref{F:Multiplication}) Using the formula~$\hh(ab)=\hh(\hh(a)b)$ from
Lemma~\ref{L:KeyLemma}, we obtain
\begin{eqnarray*}
   \hh(y_jh_j\cdots h_i)&=&\hh(\hh(y_jh_j)h_{j-1}\cdots h_i)\\
    &=&\hh(y_{j-1}h_{j-1}\cdots h_i)=\ldots=\hh(y_ih_i)=y_{i-1}
\end{eqnarray*}
for~$1\leq i<j\leq m$, and
\[
\hh(h_m\cdots h_iz_{i-1})=\hh(\hh(h_m\cdots h_i)z_{i-1})=\hh(h_iz_{i-1})=h''_i
\]
for~$1\leq i\leq m$.

The only point remaining to be checked is that the $h'_i$'s and~$z_m$ are hypercubes. For this, it suffices to show that the right normal form of the product of two hypercubes has length at most two. Assume that~$a,b$ are two hypercubes. We denote by~$a'$ the element satisfying~$ab=a'\hh(ab)$. As~$b$ is a hypercube, we have~$\hh(ab)=a''b$ for some hypercube~$a''$. By right cancellation, we obtain~$a=a'a''$. Then~$a'$ divides the hypercube~$a$, and, therefore, $a'$ is a hypercube too. This concludes the proof.\qed\end{pf}

\noindent Although quite natural, the latter result was not obvious beforehand. Indeed, putting in normal form a product of two hypercubes might have required say three hypercubes, since the condition for being normal discards some decompositions.

\noindent Proposition~\ref{P:Multiplication} is therefore exactly what we need to show that, in addition to the regularity of the language of its right normal forms, every left divisibility monoid satisfies the required deep geometric property~: 

\begin{cor}\label{C:FellowTraveller} The language of the right normal forms of every left divisibility monoid satisfies both the left and the right directed fellow traveller properties.
\end{cor}

\begin{pf} We are going to show that, for every hypercube~$y$ (\resp every hypercube~$z$), the left (\resp right) directed distance between the right normal form of an element~$a$ and the one of~$ya$ (\resp the one of~$az$) is uniformly bounded by~$2$, what will establish the left (\resp right) directed fellow traveller property.

The point is to consider the set~$\HH$ of its hypercubes as a generating alphabet for~$M$.

(i) Let~$a,b$ in~$M$ satisfy~$d_{\HH}(a,b)\leq 1$. The result is trivial for~$d_{\HH}(a,b)=0$. Assume~$d_{\HH}(a,b)=1$.  Then we can suppose, without loss of generality, that~$b=ya$ holds for some~$y$ in~$\HH$. Let~$h_m\cdot\ldots\cdot h_1$ be the right normal form of~$a$. By Proposition~\ref{P:Multiplication}(i), the right normal form of~$b$ is~$h'_m\cdot\ldots\cdot h'_1\cdot y_0$ (or possibly $h'_{m-1}\cdot\ldots\cdot h'_1\cdot y_0$) with~$y_m=y$, $y_{i-1}=\hh(y_ih_i)$ and~$y_ih_i=h'_iy_{i-1}$ for~$1\leq i\leq m$. Therefore, for every~$t>0$, we find
\[d_{\HH}(a_{[t]},b_{[t]})=d_{\HH}(h_t\cdot\ldots\cdot h_1,h'_{t-1}\cdot\ldots\cdot h'_1\cdot y_0)<2.\]
Indeed, by definition, we have~$y_{t-1}=\hh(y_th_t)$ and, as~$h_t$ is a hypercube, we have~$y_{t-1}=y'_th_t$ for some hypercube~$y'_t$~: we obtain~$y'_th_t\cdots h_1=y_{t-1}h_{t-1}\cdots h_1=h'_{t-1}\cdots h'_1y_0$.

(ii) Let~$a,b$ in~$M$ satisfy~$d^{\HH}(a,b)\leq 1$. The result is trivial for~$d^{\HH}(a,b)=0$. Assume~$d^{\HH}(a,b)=1$.  Then we can suppose, without loss of generality, that~$b=az$ holds for some~$z$ in~$\HH$. Let~$h_m\cdot\ldots\cdot h_1$ be the right normal form of~$a$. By Proposition~\ref{P:Multiplication}(ii), the right normal form of~$b$ is~$z_m\cdot h''_m\cdot\ldots\cdot h''_1$ (or possibly $h''_m\cdot\ldots\cdot h''_1$) with~$z_0=z$, $h''_i=\hh(h_iz_{i-1})$ and~$z_ih''_i=h_iz_{i-1}$ for~$1\leq i\leq m$. Two cases may occur. First, assume~$z_m=1$. Then, for every~$t>0$, we find
\[d^{\HH}(a^{[t]},b^{[t]})=d^{\HH}(h_m\cdot\ldots\cdot h_{m-t+1},h''_{m}\cdot\ldots\cdot h''_{m-t+1})<2,\]
since we have $h_m\cdots h_{m-t+1}=h''_{m}\cdots h''_{m-t+1}z_{m-t}$. Now assume~$z_m\not =1$. Then, for every~$t>0$, we find
\[d^{\HH}(a^{[t]},b^{[t]})=d^{\HH}(h_m\cdot\ldots\cdot h_{m-t+1},z_m\cdot h''_{m}\cdot\ldots\cdot h''_{m-t+2})<2.\]
Indeed, by definition, we have~$h''_{m-t+1}=\hh(h_{m-t+1}z_{m-t})$ and, as~$z_{m-t}$ is a hypercube, we have~$h''_{m-t+1}=z'_{m-t}z_{m-t}$ for some hypercube~$z'_{m-t}$~: we obtain as required~$h_m\cdots h_{m-t+1}=z_mh''_{m}\cdots h''_{m-t+2}z'_{m-t}$. This concludes the proof.\qed
\end{pf}

\section{Biautomaticity and associated finite transducers}
\label{Sec:Transducers}

\noindent In this section, we establish the main theorems of the paper and illustrate them with
several examples. We finally discuss about transducers and multiplier automata.

\subsection{Biautomaticity}

\noindent The results from the previous section make us ready to establish the biautomaticity of left divisibility monoids.

\begin{prop}\label{P:Automatic} The language of the right normal forms provides a biautomatic
structure to every left divisibility monoid.
\end{prop}

\begin{pf} According to~Theorem~\ref{P:HC}, in a cancellative monoid, every regular
language satisfying the directed fellow traveller property for both left and right multiplication
provides a biautomatic structure.\qed
\end{pf}

\noindent The latter provides an original and complete proof to the fact that every left divisibility monoid is automatic \cite{kus}. Actually, Proposition~\ref{P:Automatic} allows to state~:

\begin{thm}\label{T:Biautomatic} Every left divisibility monoid is biautomatic.
\end{thm} 

\subsection{Finite transducers computing right normal forms}

\noindent Following Thurston's original idea concerning the automaticity of the braid groups
(see \cite{eps}), Dehornoy constructed in~\cite{dhh} an explicit finite transducer computing normal forms in every Garside monoid (see also~\cite{pie}). We show here that these methods can be adapted to left divisibility monoids. All these transducers work similarly to the one in~Example~\ref{E:Base}.

\begin{thm}\label{T:Transducer} Every left divisibility monoid admits an explicit finite
transducer which allows to compute right normal forms in quadratic time.
\end{thm}
 
\begin{pf} Let~$M$ be a left divisibility monoid and~$\Irr$ the set of its irreducible elements. The
transducer can be built as follows. First, the set of the states is exactly the set~$\mathcal H$ of
the hypercubes in~$M$. Next, for every state~$a$ in~$\mathcal H$ and every irreducible~$x$
in~$\Irr$, there is an arrow from~$a$ to the state~$b$ defined to be the maximal hypercube in~$\mathcal H$
right-dividing~$ax$, which is well-defined according to Lemma~\ref{L:RightDividing}~; this arrow is
then labelled by~$x|u$ where~$u$ is any word over~$\Irr$ satisfying~$a\overline{x}=\overline{u}b$.

\noindent The right normal form~$\NN(w)$ of a given word~$w$ over~$\Irr$---formally defined to be~$\NN(\overline{w})$ according to Definition~\ref{D:RightNF}---is then computed as follows. During the reading of~$w$ by the just defined transducer, one concatenates the corresponding outputs (eventually empty, namely~$\varepsilon$) of trodden arrows. At the end of the reading of~$w$, the ambient state~$s$ is the first hypercube of the right normal form~$\NN(w)$ of~$w$ and the word~$w'$ obtained by concatenating the various outputs is the word that remains to be normalized~: we have~$\NN(w)=\NN(w')\cdot s$. With~${\bm\lambda}$ and~${\bm\tau}$ as in~Definitions~\ref{D:Auto} and~\ref{D:Trans}, we obtain~$\NN(\varepsilon)={\mathbf 1}$ and~$\NN(w)=\NN({\bm\lambda}(w,{\mathbf 1}))\cdot{\bm\tau}(w,{\mathbf 1})$ for~$w\not=\varepsilon$.\qed\end{pf}

\begin{rem} Additional arrows allow the just defined transducers to become able to read words over the whole alphabet~$\mathcal H$. Precisely, the associated \emph{augmented} transducer can be defined as follows. Again, the set of the states is~$\mathcal H$. Then, for every state~$a$ in~$\mathcal H$ and every hypercube~$h$ in~$\mathcal H$, there is an arrow from~$a$ to the state~$b$ defined to be the maximal hypercube in~$\mathcal H$ right-dividing~$ah$, which is well-defined according to Lemma~\ref{L:RightDividing}~; this arrow is then labelled by~$h|k$ where~$k$ is the hypercube satisfying~$ah=kb$, which is well-defined according to Proposition~\ref{P:Multiplication}.
\end{rem}

\subsection{Three examples}

\noindent We apply Theorem~\ref{T:Transducer} on those left divisibility monoids from Example~\ref{E:DiviMonoid} and from Remark~\ref{R:Connectivity}.

\begin{figure}[thb]
	\begin{center}
		\includegraphics[width=220pt]{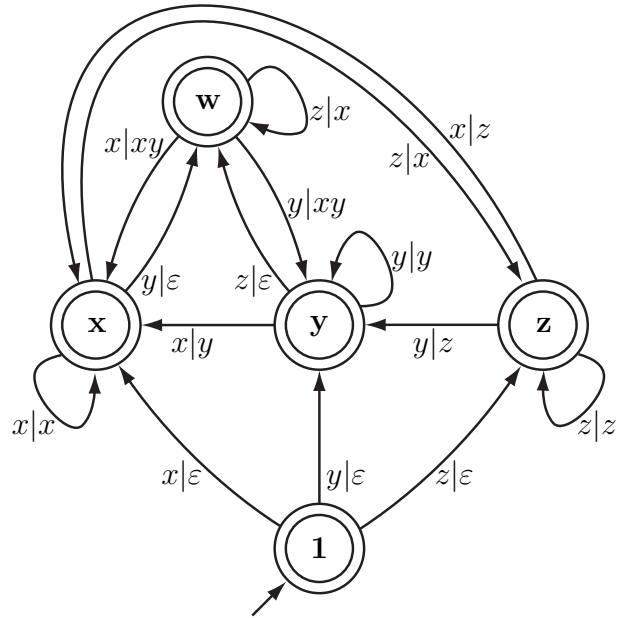}
		\put(0,245){\makebox(0,0){$$}}
		\put(-110,26){\makebox(0,0){$\mathbf{1}$}}
		\put(-194,110.5){\makebox(0,0){$\mathbf{x}$}}
		\put(-110,110){\makebox(0,0){$\mathbf{y}$}}
		\put(-25,110.5){\makebox(0,0){$\mathbf{z}$}}
		\put(-158,103){\makebox(0,0){$x|y$}}
		\put(-67,103){\makebox(0,0){$y|z$}}
		\put(-218,72){\makebox(0,0){$x|x$}}
		\put(-75,136){\makebox(0,0){$y|y$}}
		\put(-5,72){\makebox(0,0){$z|z$}}
		\put(-162,53){\makebox(0,0){$x|\varepsilon$}}
		\put(-100,53){\makebox(0,0){$y|\varepsilon$}}
		\put(-59,53){\makebox(0,0){$z|\varepsilon$}}
		\put(-170,127){\makebox(0,0){$y|\varepsilon$}}
		\put(-111,156){\makebox(0,0){$y|xy$}}
		\put(-180,179){\makebox(0,0){$x|xy$}}
		\put(-135,127){\makebox(0,0){$z|\varepsilon$}}
		\put(-152,195){\makebox(0,0){$\mathbf{w}$}}
		\put(-106,191){\makebox(0,0){$z|x$}}
		\put(-76,173){\makebox(0,0){$z|x$}}
		\put(-53,184){\makebox(0,0){$x|z$}}
	\end{center}
\caption{The transducer for~$\langle~x,y,z:xy=yz~\rangle$.}
 \label{F:Transducer1}
\end{figure}

\begin{exmp} Figure~\ref{F:Transducer1} displays the 5-state transducer which allows to compute right normal forms in the left divisibility monoid~$\langle~x,y,z:xy=yz~\rangle$. Indeed, the hypercubes are~$\mathbf{1}=\jjR\{~\}$, $\mathbf{x}=\jjR\{x\}$, $\mathbf{y}=\jjR\{y\}$, $\mathbf{z}=\jjR\{z\}$ and~$\mathbf{w}=\overline{xy}=\overline{yz}=\jjR\{x,y\}$. The initial state is the state~$\mathbf{1}$. Next, the~$5\times 3=15$ labelled arrows are constructed as in the previous proof. For instance, there is a arrow between the state~$\mathbf{w}$ and the maximal hypercube right-dividing $\mathbf{w}\overline{z}=\overline{xyz}=\overline{x}\mathbf{w}$, namely the state~$\mathbf{w}$ itself~; this arrow is then labelled by~$z|x$.
\end{exmp}

\begin{figure}[thb]
	\begin{center}
		\includegraphics[width=220pt]{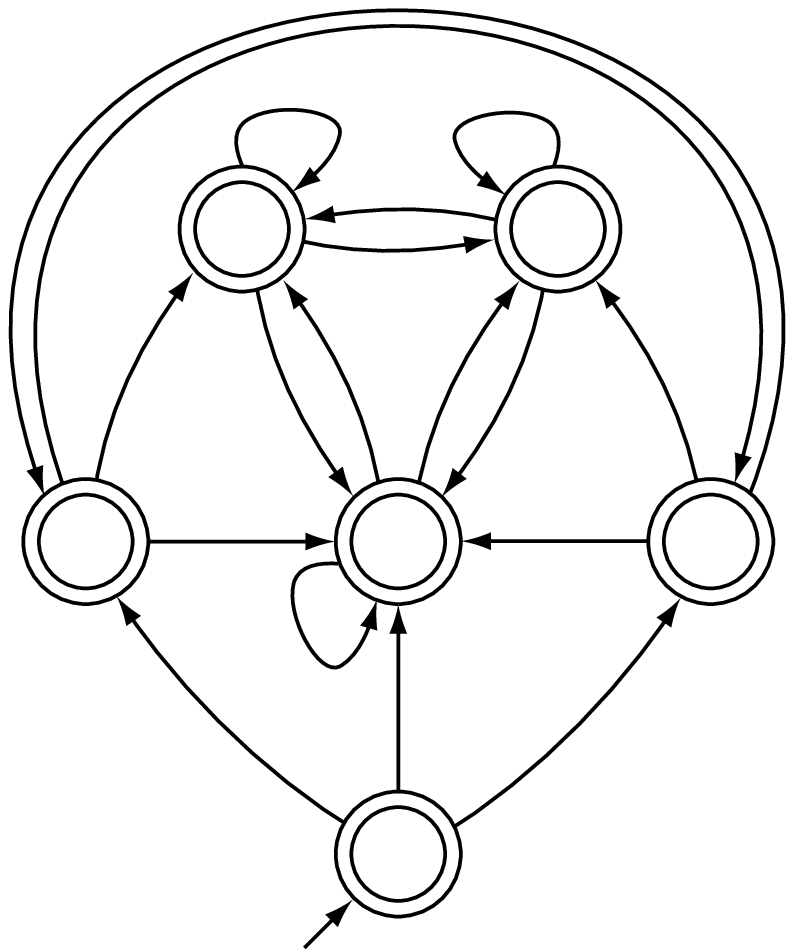}
		\put(0,270){\makebox(0,0){$$}}
		\put(-110,26){\makebox(0,0){$\mathbf{1}$}}
		\put(-194,110){\makebox(0,0){$\mathbf{x}$}}
		\put(-109,110){\makebox(0,0){$\mathbf{y}$}}
		\put(-25,110){\makebox(0,0){$\mathbf{z}$}}
		\put(-155,101){\makebox(0,0){$y|x$}}
		\put(-64,101){\makebox(0,0){$y|z$}}
		\put(-151,196){\makebox(0,0){$\mathbf{x^2}$}}
		\put(-66,196){\makebox(0,0){$\mathbf{z^2}$}}
		\put(-137,72){\makebox(0,0){$y|y$}}
		\put(-164,52){\makebox(0,0){$x|\varepsilon$}}
		\put(-100,52){\makebox(0,0){$y|\varepsilon$}}
		\put(-57,52){\makebox(0,0){$z|\varepsilon$}}
		\put(-179,133){\makebox(0,0){$x|\varepsilon$}}
		\put(-120,168){\makebox(0,0){$z|\varepsilon$}}
		\put(-99,168){\makebox(0,0){$x|\varepsilon$}}
		\put(-42,133){\makebox(0,0){$z|\varepsilon$}}
		\put(-198,171){\makebox(0,0){$z|x$}}
		\put(-225,155){\makebox(0,0){$x|z$}}
		\put(-145,133){\makebox(0,0){$y|x^2$}}
		\put(-77,133){\makebox(0,0){$y|z^2$}}
		\put(-127,232){\makebox(0,0){$x|x$}}
		\put(-111,182){\makebox(0,0){$z|y$}}
		\put(-111,208){\makebox(0,0){$x|y$}}
		\put(-92,232){\makebox(0,0){$z|z$}}
	\end{center}
\caption{The transducer for~$\langle~x,y,z:x^2=yz,yx=z^2~\rangle$.}
\label{F:Transducer2}
\end{figure}

\begin{exmp} Figure~\ref{F:Transducer2} displays the 6-state transducer for the second left divisibility monoid in Example~\ref{E:DiviMonoid}, namely the monoid~$\langle~x,y,z:x^2=yz,yx=z^2~\rangle$. Let us observe how this transducer allows us to compute the right normal form of a word---say the word~$w_0=yzyxxz$---in the associated monoid. The reading of~$w_0$ from the state~$\mathbf{1}$ leads to the state~$s_1=\mathbf{z^2}$, and the concatenation of the corresponding outputs is the word~$w_1=xxyy$. Therefore, we have~$\NN(w_0)=\NN(w_1)\cdot s_1$. The word~$w_1$ obtained is the word that remains to be normalized~: the reading of the word~$w_1=xxyy$ from the state~$\mathbf{1}$ leads to the state~$s_2=\mathbf{y}$, and the concatenation of the output labels of trodden arrows is the word~$w_2=xxy$. We obtain~$\NN(w_0)=\NN(w_2)\cdot s_2\cdot s_1$. Repeating the process twice again, we finally obtain the right normal form~$\NN(yzyxxz)=\mathbf{x^2}\cdot\mathbf{y}\cdot\mathbf{y}\cdot\mathbf{z^2}$.
\end{exmp}

\begin{figure}[thb]
	\begin{center}
		\includegraphics[width=220pt]{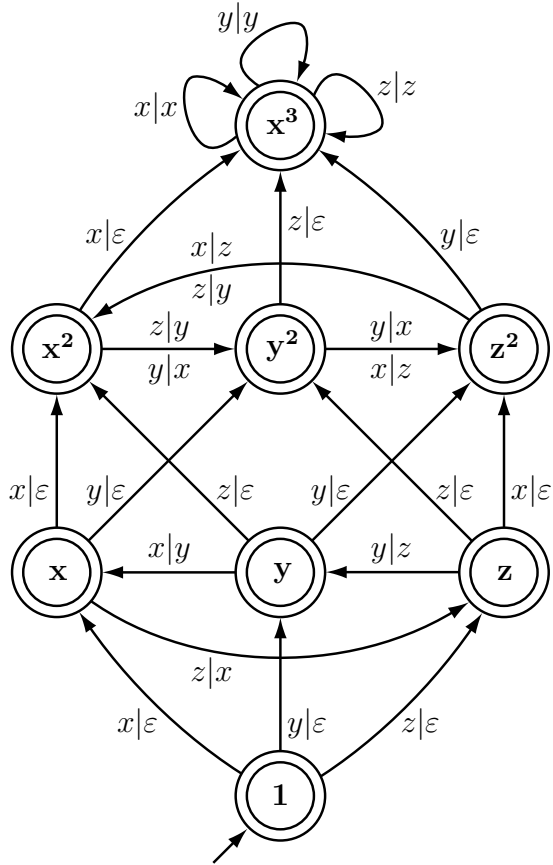}
		\put(0,330){\makebox(0,0){$$}}
		\put(-110,26){\makebox(0,0){$\mathbf{1}$}}
		\put(-194,110){\makebox(0,0){$\mathbf{x}$}}
		\put(-109,110){\makebox(0,0){$\mathbf{y}$}}
		\put(-25,110){\makebox(0,0){$\mathbf{z}$}}
		\put(-194,196){\makebox(0,0){$\mathbf{x^2}$}}
		\put(-109,196){\makebox(0,0){$\mathbf{y^2}$}}
		\put(-25,196){\makebox(0,0){$\mathbf{z^2}$}}
		\put(-109,281){\makebox(0,0){$\mathbf{x^3}$}}
		\put(-164,53){\makebox(0,0){$x|\varepsilon$}}
		\put(-100,53){\makebox(0,0){$y|\varepsilon$}}
		\put(-57,53){\makebox(0,0){$z|\varepsilon$}}
		\put(-205,141){\makebox(0,0){$x|\varepsilon$}}
		\put(-176,141){\makebox(0,0){$y|\varepsilon$}}
		\put(-127,141){\makebox(0,0){$z|\varepsilon$}}
		\put(-91,141){\makebox(0,0){$y|\varepsilon$}}
		\put(-44,141){\makebox(0,0){$z|\varepsilon$}}
		\put(-15,141){\makebox(0,0){$x|\varepsilon$}}
		\put(-176,238){\makebox(0,0){$x|\varepsilon$}}
		\put(-100,243){\makebox(0,0){$z|\varepsilon$}}
		\put(-42,238){\makebox(0,0){$y|\varepsilon$}}
		\put(-152,119){\makebox(0,0){$x|y$}}
		\put(-136,72){\makebox(0,0){$z|x$}}
		\put(-68,119){\makebox(0,0){$y|z$}}
		\put(-152,203){\makebox(0,0){$z|y$}}
		\put(-152,187){\makebox(0,0){$y|x$}}
		\put(-68,203){\makebox(0,0){$y|x$}}
		\put(-68,187){\makebox(0,0){$x|z$}}
		\put(-136,233){\makebox(0,0){$x|z$}}
		\put(-136,217){\makebox(0,0){$z|y$}}
		\put(-66,294){\makebox(0,0){$z|z$}}
		\put(-126,320){\makebox(0,0){$y|y$}}
		\put(-157,286){\makebox(0,0){$x|x$}}
	\end{center}
\caption{The transducer for~$\langle~x,y,z:x^2=yz,y^2=zx,z^2=xy~\rangle$.}
\label{F:Transducer3}
\end{figure}

\begin{exmp} Figure~\ref{F:Transducer3} displays the 8-state transducer allowing to compute right normal forms in the left divisibility monoid $\langle~\xx,\yy,\zz:\xx^2=\yy\zz,\yy^2=\zz\xx,\zz^2=\xx\yy~\rangle$ introduced in Remark~\ref{R:Connectivity}. The transducer can be compared with the corresponding graph of hypercubes of Figure~\ref{F:Graph}. Moreover, this example of a divisibility monoid distinguishes from each of the two latter by its non width-boundedness, and so by its non rationality (see~\cite{kus}). Indeed, for every~$n\in\mathbb{N}$, the width of the lattice~$\Dar(\xx^n)$ equals the number of partitions of~$n$ into at most three parts, hence equals the nearest integer to~$\frac{(n+3)^2}{12}$ (sequence~{\rm A001399} from~\cite{slo}).
\end{exmp}

\begin{figure}[thb]
	\begin{center}
		\includegraphics[width=220pt]{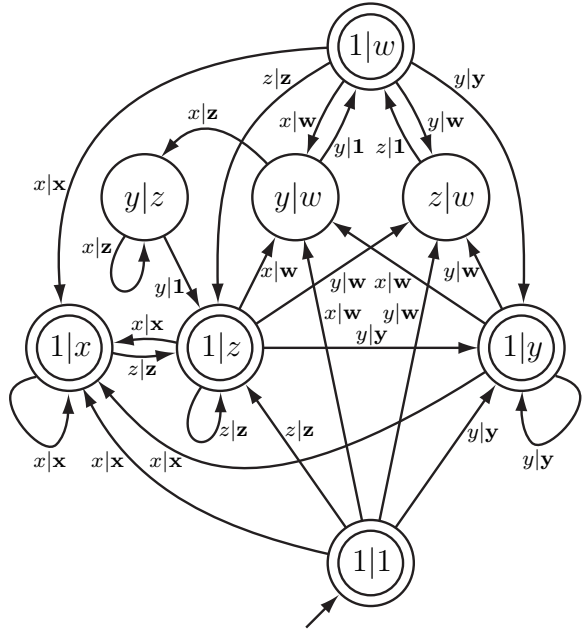}
		\put(-81,25){\makebox(0,0){$1|1$}}
		\put(-195,106){\makebox(0,0){$1|x$}}
		\put(-138,106){\makebox(0,0){$1|z$}}
		\put(-24,106){\makebox(0,0){$1|y$}}
		\put(-167,163){\makebox(0,0){$y|z$}}
		\put(-109,163){\makebox(0,0){$y|w$}}
		\put(-51,163){\makebox(0,0){$z|w$}}
		\put(-81,220){\makebox(0,0){$1|w$}}
		\put(-203,63){\makebox(0,0){\scriptsize$x|{\mathbf x}$}}
		\put(-181,63){\makebox(0,0){\scriptsize$x|{\mathbf x}$}}
		\put(-159,63){\makebox(0,0){\scriptsize$x|{\mathbf x}$}}
		\put(-166,115){\makebox(0,0){\scriptsize$x|{\mathbf x}$}}
		\put(-203,168){\makebox(0,0){\scriptsize$x|{\mathbf x}$}}
		\put(-168,98){\makebox(0,0){\scriptsize$z|{\mathbf z}$}}
		\put(-133.5,75){\makebox(0,0){\scriptsize$z|{\mathbf z}$}}
		\put(-108.5,75){\makebox(0,0){\scriptsize$z|{\mathbf z}$}}
		\put(-157.5,129){\makebox(0,0){\scriptsize$y|{\mathbf 1}$}}
		\put(-185,144){\makebox(0,0){\scriptsize$x|{\mathbf z}$}}
		\put(-39,74){\makebox(0,0){\scriptsize$y|{\mathbf y}$}}
		\put(-18,64){\makebox(0,0){\scriptsize$y|{\mathbf y}$}}
		\put(-44,208){\makebox(0,0){\scriptsize$y|{\mathbf y}$}}
		\put(-117,208){\makebox(0,0){\scriptsize$z|{\mathbf z}$}}
		\put(-81,111){\makebox(0,0){\scriptsize$y|{\mathbf y}$}}
		\put(-92,119){\makebox(0,0){\scriptsize$x|{\mathbf w}$}}
		\put(-70,119){\makebox(0,0){\scriptsize$y|{\mathbf w}$}}
		\put(-73,131){\makebox(0,0){\scriptsize$x|{\mathbf w}$}}
		\put(-90,131){\makebox(0,0){\scriptsize$y|{\mathbf w}$}}
		\put(-116,135){\makebox(0,0){\scriptsize$x|{\mathbf w}$}}
		\put(-47,135){\makebox(0,0){\scriptsize$y|{\mathbf w}$}}
		\put(-89.5,182){\makebox(0,0){\scriptsize$y|{\mathbf 1}$}}
		\put(-74,182){\makebox(0,0){\scriptsize$z|{\mathbf 1}$}}
		\put(-109,192){\makebox(0,0){\scriptsize$x|{\mathbf w}$}}
		\put(-53,192){\makebox(0,0){\scriptsize$y|{\mathbf w}$}}
		\put(-145,195){\makebox(0,0){\scriptsize$x|{\mathbf z}$}}
	\end{center}
\caption{Kuske's transducer for~$\langle~x,y,z:xy=yz~\rangle$.}
 \label{F:TransducerK}
\end{figure}

\begin{rem}\label{R:Kuske} The transducer of Figure~\ref{F:Transducer1} can be compared to Kuske's transducer of Figure~\ref{F:TransducerK} whose relevance was mentioned at the very end of Section~\ref{Sec:Background2}. Recall that, in this rational case, a normal form is computed through only one run.
\end{rem}

\subsection{Transducers \emph{vs} multiplier automata}\label{SS:vs}

\noindent In the context of computing normal forms in divisibility monoids, the afore described transduction machinery seems to afford several advantages over the classical multiplier automata. These advantages revolve around two main ideas~: legibility and efficiency.
The legibility of our transducer comes from the relative compactness of the data, and, above all, from the fact that its graph structure mimics the lattice structure that the set of hypercubes is endowed with. 

As regards efficiency, one can observe that, even if the multiplier automata---which are chosen deterministic---can be viewed and used as transducers, the latter are in general neither subsequential nor even subsequentiable. The time efficiency is known to be substantially increased when subsequential machines are used (see~\cite{moh} for instance). Moreover, it is worth to noticing that our transducers are able to compute the rightmost hypercube of the right normal form in linear time.

\begin{figure}[thb]
	\begin{center}
		\includegraphics[width=400pt]{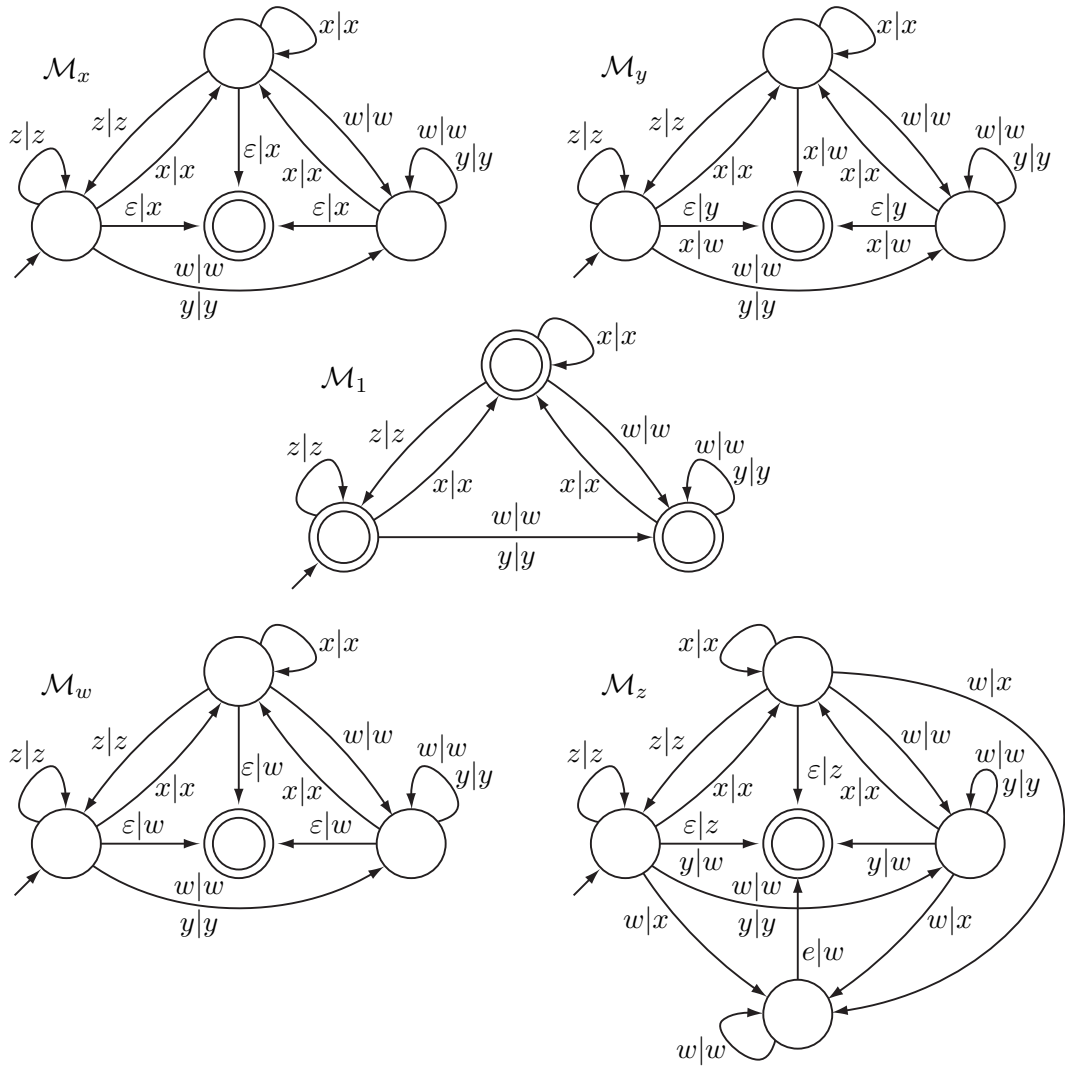}
		\put(-378,378){\makebox(0,0){\small $\mathcal{M}_x$}}
		\put(-393,353){\makebox(0,0){\small $z|z$}}
		\put(-362,358){\makebox(0,0){\small $z|z$}}
		\put(-337,340){\makebox(0,0){\small $x|x$}}
		\put(-349,326){\makebox(0,0){\small $\e|x$}}
		\put(-304,348){\makebox(0,0){\small $\e|x$}}
		\put(-289,339){\makebox(0,0){\small $x|x$}}
		\put(-278,326){\makebox(0,0){\small $\e|x$}}
		\put(-275,395){\makebox(0,0){\small $x|x$}}
		\put(-264,360){\makebox(0,0){\small $w|w$}}
		\put(-236,355){\makebox(0,0){\small $w|w$}}
		\put(-224,345){\makebox(0,0){\small $y|y$}}
		\put(-328,303){\makebox(0,0){\small $w|w$}}
		\put(-328,289){\makebox(0,0){\small $y|y$}}
		\put(-273,261){\makebox(0,0){\small $\mathcal{M}_1$}}
		\put(-288,235){\makebox(0,0){\small $z|z$}}
		\put(-257,240){\makebox(0,0){\small $z|z$}}
		\put(-232,221){\makebox(0,0){\small $x|x$}}
		\put(-208,209){\makebox(0,0){\small $w|w$}}
			\put(-208,194){\makebox(0,0){\small $y|y$}}
		\put(-184,221){\makebox(0,0){\small $x|x$}}
		\put(-170,277){\makebox(0,0){\small $x|x$}}
		\put(-159,242){\makebox(0,0){\small $w|w$}}
		\put(-131,236){\makebox(0,0){\small $w|w$}}
		\put(-119,226){\makebox(0,0){\small $y|y$}}
		\put(-167,378){\makebox(0,0){\small $\mathcal{M}_y$}}
		\put(-182,353){\makebox(0,0){\small $z|z$}}
		\put(-151,358){\makebox(0,0){\small $z|z$}}
		\put(-126,340){\makebox(0,0){\small $x|x$}}
		\put(-138,326){\makebox(0,0){\small $\e|y$}}
			\put(-138,312){\makebox(0,0){\small $x|w$}}
		\put(-91,347){\makebox(0,0){\small $x|w$}}
		\put(-78,339){\makebox(0,0){\small $x|x$}}
		\put(-67,326){\makebox(0,0){\small $\e|y$}}
			\put(-67,312){\makebox(0,0){\small $x|w$}}
		\put(-64,395){\makebox(0,0){\small $x|x$}}
		\put(-53,360){\makebox(0,0){\small $w|w$}}
		\put(-25,355){\makebox(0,0){\small $w|w$}}
		\put(-13,345){\makebox(0,0){\small $y|y$}}
		\put(-117,303){\makebox(0,0){\small $w|w$}}
		\put(-117,289){\makebox(0,0){\small $y|y$}}
		\put(-378,145){\makebox(0,0){\small $\mathcal{M}_w$}}
		\put(-393,119){\makebox(0,0){\small $z|z$}}
		\put(-362,124){\makebox(0,0){\small $z|z$}}
		\put(-337,106){\makebox(0,0){\small $x|x$}}
		\put(-349,92){\makebox(0,0){\small $\e|w$}}
		\put(-304,114){\makebox(0,0){\small $\e|w$}}
		\put(-289,105){\makebox(0,0){\small $x|x$}}
		\put(-278,92){\makebox(0,0){\small $\e|w$}}
		\put(-275,161){\makebox(0,0){\small $x|x$}}
		\put(-264,126){\makebox(0,0){\small $w|w$}}
		\put(-236,121){\makebox(0,0){\small $w|w$}}
		\put(-224,111){\makebox(0,0){\small $y|y$}}
		\put(-328,69){\makebox(0,0){\small $w|w$}}
		\put(-328,55){\makebox(0,0){\small $y|y$}}
		\put(-167,145){\makebox(0,0){\small $\mathcal{M}_z$}}
		\put(-182,119){\makebox(0,0){\small $z|z$}}
		\put(-151,124){\makebox(0,0){\small $z|z$}}
		\put(-126,106){\makebox(0,0){\small $x|x$}}
		\put(-138,92){\makebox(0,0){\small $\e|z$}}
			\put(-138,78){\makebox(0,0){\small $y|w$}}
		\put(-91,113){\makebox(0,0){\small $\e|z$}}
		\put(-78,105){\makebox(0,0){\small $x|x$}}
			\put(-67,78){\makebox(0,0){\small $y|w$}}
		\put(-139,161){\makebox(0,0){\small $x|x$}}
		\put(-53,126){\makebox(0,0){\small $w|w$}}
			\put(-29,146){\makebox(0,0){\small $w|x$}}
		\put(-25,119){\makebox(0,0){\small $w|w$}}
		\put(-17,108){\makebox(0,0){\small $y|y$}}
		\put(-117,69){\makebox(0,0){\small $w|w$}}
		\put(-117,55){\makebox(0,0){\small $y|y$}}
			\put(-159,56){\makebox(0,0){\small $w|x$}}
			\put(-92,44){\makebox(0,0){\small $e|w$}}
			\put(-44,56){\makebox(0,0){\small $w|x$}}
			\put(-138,7){\makebox(0,0){\small $w|w$}}
	\end{center}
\caption{Multiplier automata for~$\langle~x,y,z:xy=yz~\rangle$.}
\label{F:Multipliers}
\end{figure}

\noindent In order to illustrate the purpose, the reader is invited to compare Figures~\ref{F:Transducer1} and~\ref{F:Multipliers}. The latter shows the equality recognizer automaton~$\mathcal{M}_1$ and the multiplier automata~$ \mathcal{M}_\xx$, $\mathcal{M}_\yy$, $\mathcal{M}_\zz$ and~$\mathcal{M}_\ww$. In order to compute the right normal form of a word~$x_1\cdots x_n$ over~$\Irr$ or even over~$\HH$, the right normal forms of successive prefixes~$x_1\cdots x_i$ are computed by applying~$\mathcal{M}_{x_i}$ to the normal form~$x_1\cdots x_{i-1}$.

\noindent Finally, it appears that, in the case of Artin's braid monoids, of Garside monoids, of free partially commutative monoids, and now in the case of divisibility monoids, the biautomaticity and the associated transducers are designed from the particular structure---namely a lattice or semi-lattice structure---of the divisibility relation. Even if a global approach seems to be out of reach, few experimental investigations indicate that wider classes of automatics monoids could be study with similar tools and afford a new insight on the subject.

\begin{ack}
The author is grateful to Manfred Droste, Dietrich Kuske, Sylvain Lombardy, Jean Mairesse and Jacques Sakarovitch for several helpful conversations and the anonymous referees for their careful readings of the paper.
\end{ack}


\end{document}